\documentclass[a4paper,11pt]{article} 
 
\usepackage{fancybox} 
\usepackage{pifont} 

\usepackage[english]{babel}

\usepackage[symbol]{footmisc}

\usepackage{dsfont}
\usepackage{amsmath} 
\usepackage[applemac]{inputenc}
\usepackage{amsfonts}
\usepackage{amssymb}
\usepackage{amsthm}
\usepackage{stmaryrd}
\usepackage{color}
\usepackage{mathrsfs}

\usepackage[backref]{hyperref}

\newcommand{\mysection}{\setcounter{equation}{0} \section}

\renewcommand{\P}{\mathbb{P}}

\newcommand{\N}{\mathbb{N}}

\newcommand{\R}{\mathbb{R}}

\newtheorem{THM}{Theorem}   
\newtheorem{remark}{Remark}

\def\1{\mbox{1\hspace{-0.25em}l}}

\def\leftB{[\![}
\def\rightB{]\!]}

\def\det{{{\rm det}}}

\def\0{{\mathbf{0}}}

\newcommand{\newabstract}[1]{%
	\par\bigskip
	\csname otherlanguage*\endcsname{#1}%
	\csname captions#1\endcsname
	\item[\hskip\labelsep\scshape\abstractname.]
}

\title{\textbf{%Régularisation par changement d'échelle temporelle
	When Pythagoras meets Navier-Stokes}}
 \author{\textbf{Igor Honor\'e}\footnote{Univ Lyon, CNRS, Universit\'e Claude Bernard Lyon 1, UMR5208, Institut Camille Jordan, F-69622 Villeurbanne, France. E-mail: honore@math.univ-lyon1.fr}}

\begin{document}

\maketitle

%\begin{abstract}

% \selectlanguage{english}

 \begin{abstract}
 In this article, we develop a new method,  based on a time decomposition of a Cauchy problem elaborated in  \cite{hono_transport_v2:22}, to retrieve the well-known  $L^\infty ([0,T],L^2(\R^d,\R^d))$ control of the solution of the incompressible Navier-Stokes equation in $\R^d$. We precisely explain how the Pythagorean theorem in $L^2(\R^d,\R^d)$ allows to get the proper energy estimate; however such an argument does not work anymore in $L^p(\R^d,\R^d)$, $p \neq 2$.
 We also deduce, by similar arguments, an already known $L^\infty ([0,T],L^1(\R^3,\R^3))$ control of vorticity for $d=3$.
% ,
% % $u$ of the incompressible Navier-Stokes equation in $\R^3$ 
% already proved by \cite{cons:90}.
% \textcolor{red}{Generaliser à une equation plus generale avec des projections ? Pourquoi pas dans un corollaire avec supposition de l'existence, ou pas ... EN tous cas bien preciser dans les comms}
% Our method does rely on the Leray's approach, but on the time decomposition technique introduced in \cite{hono_transport_v2:22}.
 
 \end{abstract}
 
 {\small{\textbf{Keywords:} Incompressible Navier-Stokes equation, Energy estimates.}}
 
 {\small{\textbf{MSC:} Primary: 35K40, 35Q35; Secondary: 46E35, 76D05 }}

 \mysection{Introduction}

For each $d \in \N$, we consider the incompressible Navier-Stokes equation  defined by
 \begin{equation}
 \label{Navier_Stokes_equation_v1}
 \begin{cases}
 \partial_t  u(t,x)+  u(t,x)   \cdot  \nabla   u(t,x)
 = \nu \Delta u(t,x)-\nabla {\rm p}(t,x)+  f(t,x) 
 ,\\
 u(0,x)= u_0(x), \\
 \nabla \cdot  u(t,x)= 0, 
 \end{cases}
 \end{equation}
 where $ f, u_0$ satisfy some regularity conditions.
 % of \cite{osee:11}, \cite{osee:12}.
 %, or   Leray \cite{lera:34}.
 Projecting on the space of divergence free functions equivalently gives 
 \begin{equation}
 \label{Navier_Stokes_equation_v2}
 \begin{cases}
 \partial_t  u(t,x)+ \mathbb P[ u   \cdot  \nabla   u](t,x)
 = \nu \Delta u(t,x)+ \mathbb P  f(t,x) 
 ,\\
 u(0,x)= u_0(x), \, x \in \R^d ,
 \end{cases}
 \end{equation}
with $\P= id-\Xi$, the Leray operator, %which is defined by
%$\P= id-\Xi$ , 
$I $ stands for the identity operator and $\Xi$  is a non-local operator deriving from a gradient.
 
% The seminal work in small time is Oseen \cite{osee:11}, \cite{osee:12} when the time horizon is small enough.
 In three dimension, Leray \cite{lera:34} was the first to establish existence of a solution satisfying energy estimate. In this article, we use the time decomposition trick of \cite{hono_transport_v2:22} associated with a Duhamel like representation used in \cite{hono:NS_24}, to recover this control.

 The paper is organised as follows.
 We first introduce some notations and some properties of the heat kernel in Section \ref{sec_def}.
 The first result and its proof are stated in  Section \ref{sec_ain}.
% The analysis is developed in  Section \ref{sec_preuve}. 
 We dedicate Section \ref{sec_comm_pyth}  to comments on the proof about the limitation to the $L^2$ space, where Pythagorean theorem can be used, in our method.
Eventually, in Section \ref{theo_vorticity}, we adapt the analysis to prove that $\|\nabla \times u \|_{L^\infty_TL^1} < +\infty$ in $\R^3$.
% 
% \textcolor{red}{Refaire avec la norme Besov module de version continuité $B_{2,\infty}^{\gamma}$ $ \gamma \in (0,1)$}
% 
\mysection{Notations and definitions}
\label{sec_def}

We denote by $C>1$, depending only on the dimension $d$,
% (or possibly on a specific parameter $\varepsilon$ specified latter), 
a generic constant which can be different from line to line.

\subsection{Tensor and differential notations}

We use the usual differential notations, in particular, we have,
%
%The derivative operator  $\nabla$  matches with the gradient, a
%
% or with the Jacobian matrix.
%The symbol  $\nabla \cdot$ stands for the divergence,
%the Laplacian is as usual noted by $\Delta$.
%
%From these notations, we have
%Précisons maintenant les notations utilisées dans \eqref{Navier_Stokes_equation_v2}.
%Pour tout 
for any  $(t,x)\in \R_+\times \R^{d}$, 
\begin{equation*}%\label{def_contract_tensor}
u(t,x)\cdot  \nabla u(t,x) = \Big ( \sum_{j=1}^3 u_j (t,x)  \partial_{x_j}  u_i (t,x) \Big )_{i \in   \leftB 1,3 \rightB}, 
\end{equation*}
for %$b= \big ((b)_i \big)_{i \in  \leftB 1,r \rightB}$, 
$  u= \big (u_i \big)_{i \in  \leftB 1,3 \rightB}$. The Leray projector $\mathbb P$ is defined for any function $\varphi : \R^d \rightarrow \R^d$, integrable enough, by
\begin{equation}\label{def_P}
\forall x \in \R^d,	\P \varphi (x)= \varphi(x)+ \nabla (-\Delta)^{-1}\nabla\cdot \varphi(x)=:\varphi(x)-\Xi \varphi(x).
\end{equation}
In case of Navier-Stokes equation, the operator $\Xi$ matches with the pressure part 
which is indeed a gradient term in the  Helmholtz-Hodge decomposition.

\subsection{Heat kernel}
\label{sec_Gaussian_properties}

The heat kernel is denoted  for any $(t,x) \in \R_+ \times \R^d$, by
\begin{equation}\label{def_h_tx}
h_\nu (t,x):= \frac{e^{-\frac{|x|^2}{4 \nu t }}}{(4 \pi \nu t)^{\frac{3}{2}}},
\end{equation}
which satisfies  to the PDE
\begin{equation*}%\label{eq_chaleur}
\begin{cases}
(\partial_t- \nu \Delta) h_\nu(t,\cdot)=0, \forall t >0,
\\
h_\nu(0,\cdot )= \delta_0,
\end{cases}
\end{equation*}
where $\delta_0$ is the Dirac distribution.
%Pour simplifier, % les notations, 
We also denote for all $0 \leq s \leq t \leq T$ and $(x,y ) \in \R^d$
\begin{equation*}%\label{def_tilde_p}
\tilde p (s,t,x,y) = h_\nu (t-s,x-y).
\end{equation*}
The associated semi-group $\tilde P$ and  Green operator $ \tilde G$ with the heat equation are defined, for all $\varphi : \R^d \to \R^d$ and $\psi : [0,T] \times \R^d \to \R^d $ smooth enough, for any $(r,t,x) \in [0,T]^2 \times \R^d$, by :
\begin{eqnarray}\label{def_tilde_G_tilde_P}
	\tilde P_r \varphi(t,x) &:=& \int_{\R^d} \tilde p (s,t,x,y) \varphi (y) dy,
	\nonumber \\
		\tilde G_r \psi(t,x) &:=& \int_r^t  \int_{\R^d} \tilde p (s,t,x,y) \psi  (s,y) dy \, ds.
\end{eqnarray}
Furthermore, we will thoroughly use the exponential absorbing inequality: for any $\delta >0$, there a constant $ C_\delta >1$
% C \leq c>1$ 
s.t.
\begin{equation}\label{ineq_absorb}
\forall x \in \R^d, \ |x|^\delta e^{-|x|^2} \leq C_\delta e^{-C_\delta^{-1} |x|^2}.
\end{equation}

\mysection{Energy estimate}
\label{sec_ain}

%\subsection{Statement}

\begin{THM}\label{THEO}
	
%	\textcolor{red}{Dire que l'on suppose une unique soution rguliere, typiquement condition d'Oseen, ou petitesse de $u_0$.}
%	
Suppose that for a given  $T>0$, the solution $u$ of \eqref{Navier_Stokes_equation_v2} is smooth enough, %	We mean by 	smooth when 
	specifically there is $\varepsilon \in (0,1)$, such that $\|u\|_{L^\infty_T L^\infty} $, $	\|\nabla u\|_{L^\infty_T L^\infty}$, $\|\nabla  u \|_{L^\infty_T L^2}$, $
\| \nabla^2   u    \|_{L^\infty_T L^{2-\varepsilon}}  
$,
$\| \nabla	\Xi[    u       \cdot  \nabla   u   ] \|_{L^\infty_T L^\infty}$, 
$\| \nabla^2 	\Xi[    u       \cdot  \nabla   u   ] \|_{L^\infty_T L^\infty}$
$< + \infty$, then
	\begin{equation*}%\label{ineq_NS_THEO}
	\| u (T,\cdot)\|_{L^2}
	\leq 
	% e^{\frac{t}{2} }\Big (
	\|u_0\|_{L^2} +
	%\frac{1}{2}
	\int_0^T \| f(s,\cdot) \|_{L^2}  ds.
	% \Big ).
\end{equation*}

\end{THM}

\begin{remark}
When $T$ is small enough Oseen \cite{osee:11}, \cite{osee:12} shows that the solution is indeed smooth. We can iterate the above result for generic time horizon $T$.
An other possibility is to consider a ``small" $u_0$, see for instance \cite{fuji:kato:64}, \cite{giga:86}.

%It is possible to suppose less regularity on $u$, but in this case, the convergence of our approximation method would be slower.
\end{remark}

\subsection{Proof}
\label{sec_preuve}  

\subsubsection{Building of the \textit{proxy}}

%From Onseen \cite{osee:11} \cite{osee:12}, we know that there is a time horizon $T^*$ such that the solution  $u$  is smooth in $[0,T^*) \times \R^d$.

Under the hypothesis of regularity of $u$, we can set the associated ODE, defined for any \textit{freezing} point $(\tau,\xi) \in [0,T] \times  \R^d$,
\begin{equation}\label{def_theta}
	\theta_{s,\tau} (x):= x+ \int_s^\tau     u(\tilde s,\theta_{\tilde s,\tau}(x)) d \tilde s, \ s \in [0,\tau] .
\end{equation}
That is to say, for any $t \in [0,\tau]$,
\begin{equation*}%\label{def_theta}
	\dot \theta_{t,\tau} (\xi)=-    u(t,\theta_{t,\tau}(\xi)) , \ 	 \theta_{\tau ,\tau} (\xi) =\xi .
\end{equation*}
With this notation at hand, we rewrite  equation \eqref{Navier_Stokes_equation_v2} by
\begin{equation}
	\label{KOLMOLLI_xi}
	\partial_t  u (t,x)
	+    u  (t,\theta_{t,\tau} (\xi))\cdot  \nabla   u  (t,x) -   \nu \Delta   u  (t,x)
	=
	%	[ b_m(t,\xi)- b_m]
	u_{\Delta} [\tau,\xi] %(t,x) 
	\cdot  \nabla   u   (t,x) +\Xi [ u  \cdot  \nabla   u   ](t,x)
	%,
	%	+[ b_m- b] \cdot  \nabla   u (t,x)  
	+ \P   f(t,x),
\end{equation}
where
\begin{equation*}%\label{def_u_Delta}
	u_{\Delta} [\tau,\xi](t,x) := 	   u(t,\theta_{t,\tau} (\xi))  - 	   u (t,x). 
\end{equation*}
The fundamental solution of the l.h.s. in \eqref{KOLMOLLI_xi} is the perturbed heat kernel:
\begin{equation*}%\label{def_hat_p}
	\hat{p}^{\tau,\xi}  (s,t,x,y)
	:= \frac{1}{(4\pi \nu (t-s))^{\frac d 2} } 
	\exp \bigg ( -\frac {\left |x+ \int_s^t     u (\tilde s,\theta_{\tilde s,\tau }(\xi))d \tilde s-y \right|^2}{4\nu(t-s)} \bigg ). 
	%\exp\left (-\frac{|A_{s,t}(\xi)^{-1/2}(x-y)|^2}{4 }\right)
\end{equation*}
It is direct, from  definition \eqref{def_theta}, if $(\tau,\xi) =(t,x)$, that 
\begin{equation*}
	\hat p^{t,x}  (s,t,x,y)
	= \frac{1}{(4\pi \nu (t-s))^{\frac d 2} } 
	\exp \bigg ( -\frac { |\theta_{s,t } (x)-y |^2}{4\nu(t-s)} \bigg ). 
	%\exp\left (-\frac{|A_{s,t}(\xi)^{-1/2}(x-y)|^2}{4 }\right)
\end{equation*}
Also, for all 
 $\gamma \in [0,1]$, $\alpha\in \N^d$ there is a constant $C>1$ ,
\begin{eqnarray*}%\label{FIRST_deriv_CTR_DENS_flot_absorb}
		|D^\alpha \hat  p^{t,x}(s,t,x,y)|  \times  \big |y-x-\int_s^t    u (\tilde s,\theta_{\tilde s,\tau } (\xi))d \tilde s  \big |^\gamma \Big |_{(\tau,\xi)=(t,x)}
%	\nonumber \\
\hspace{-0.25cm}	&=&
	|D^\alpha \hat  p^{t,x}(s,t,x,y)| \times  \big |y-\theta_{ s,t } (x) \big |^\gamma
	%	&\leq&  C_{\alpha}
	%	\Big ( \nu (s-t)\Big ) ^{-\frac {|\alpha|}2}\bar p_{C_{\alpha}^{-1}\nu  } \Big (s,t,x+\int_0^s  _m (\tilde s,\theta_{\tilde s}^m(\xi))d \tilde s,y \Big )
	\nonumber \\
	&\leq&  C  [\nu (s-t)] ^{-\frac {|\alpha|}2+ \frac{\gamma}{2}} h_{C\nu} (t-s,\theta_{ s,t }(x)-y),%\bar p^{t,x} (s,t,x,y ),
\end{eqnarray*}
%with $\bar p^{t,x} (s,t,x,y ):=h_{C\nu} (t-s,\theta_{ s,t }(x)-y)$ 
from exponential absorbing property \eqref{ineq_absorb}.

Let us denote the Green operator for all $0 \leq r \leq t \leq  T$ and $x \in \R^d$, 
\begin{equation}\label{def_hat_G}
	\hat  G_r^{\tau,\xi}   f(t,x):= \int_r^{t}  \int_{\R^{d}}  \hat {p}^{\tau,\xi} (s,t, x,y)   f (s,y) \ dy \ ds,
\end{equation}
also the semi-group 
\begin{equation*}%\label{def_hat_P}
	\hat  P_r^{\tau,\xi}    g(t,x):=\int_{\R^{d}} \hat  p^{\tau,\xi}  (r,t,x,y)   g(y) \ dy.
\end{equation*}
We readily obtain the Duhamel formula corresponding to \eqref{KOLMOLLI_xi},
%précédente, en isolant la partie en gradient dans la décomposition d'Helmholtz-Hodge,
\begin{equation}\label{Duhamel_u_INITIAL}
	u    (t,x) 
	=  \hat P^{\tau,\xi}_{0}  u_{0}   (t,x)
	+  \hat G^{\tau,\xi}_{0} \P f(t,x) 
	+ \hat  G^{\tau,\xi}_{0}  \big (    u _{\Delta}  [\tau,\xi]  \cdot  \nabla    u    \big )(t,x)
	+ \hat  G^{\tau,\xi}_{0}  \Xi[    u      \cdot  \nabla   u   ](t,x)	
	.
\end{equation}
The last term in the r.h.s. is not \textit{a priori} divergent free when we pick $\xi=x$, (the perturbation of the heat kernel destroys the convolution structure), in order to isolate the gradient term $\Xi $, we rewrite the Duhamel formula by
\begin{eqnarray*}%\label{Duhamel_u_FINAL}
	u    (t,x) 
	&=&  \hat P^{\tau,\xi}_{0}  u_{0}   (t,x)
	+  \hat G^{\tau,\xi}_{0} \P f(t,x) 
	+ \hat  G^{\tau,\xi}_{0}  \big (    u _{\Delta}  [\tau,\xi]  \cdot  \nabla    u    \big )(t,x)
	\nonumber \\
	&&+ \hat  G^{\tau,\xi}_{0} \Big ( \Xi[    u      \cdot  \nabla   u   ](\cdot,\cdot)-\Xi[    u      \cdot  \nabla   u   ](\cdot, \theta _{\cdot,t}(x)) \Big )(t,x)	
	%	\nonumber \\
	%	&&
	+ \hat  G^{\tau,\xi}_{0} \Xi[    u      \cdot  \nabla   u   ](\cdot,\theta _{\cdot,t}(x)) (t,x)	
	\nonumber \\
	&=&  \hat P^{\tau,\xi}_{0}  u_{0}    (t,x) 
	+  \hat G^{\tau,\xi}_{0} \P f(t,x) 
	+ \hat  G^{\tau,\xi}_{0}  \big (    u_{\Delta}   [\tau,\xi]  \cdot  \nabla    u     \big )(t,x)
	\nonumber \\
	&&+ \hat  G^{\tau,\xi}_{0} \Big ( \Xi[    u       \cdot  \nabla   u    ](\cdot,\cdot)-\Xi[ u   \cdot  \nabla   u    ](\cdot,\theta _{\cdot,t}(x)) \Big )(t,x)	
	+ \int_0^t  \Xi[     u       \cdot  \nabla   u    ](s,\theta _{s,t}(x)) ds
	,
\end{eqnarray*}
where we use the obvious identity $\hat  G^{\tau,\xi}_{0} 1= \int_0^t ds$.
The last term still does not vanish after applying the Leray operator $\P$, but we can take advantage of the time integral in the flow definition.
% de divergence nulle, à cause du flot, nous ajoutons donc un terme supplémentaire.
\begin{eqnarray*}%\label{Duhamel_u_FINAL}
	u     (t,x) &=&  \hat P^{\tau,\xi}_{0}  u_{0}    (t,x) 
	+  \hat G^{\tau,\xi}_{0} \P f(t,x) 
	+ \hat  G^{\tau,\xi}_{0}  \big (    u_{\Delta}   [\tau,\xi]  \cdot  \nabla    u     \big )(t,x)
	\nonumber \\
	&&+ \hat  G^{\tau,\xi}_{0} \Big ( \Xi[    u       \cdot  \nabla   u    ](\cdot,\cdot)-\Xi[ u   \cdot  \nabla   u    ](\cdot,\theta _{\cdot,t}(x)) \Big )(t,x)	
	%\nonumber 
	\\
	&& + \int_0^t \Big ( \Xi[     u       \cdot  \nabla   u    ](s,\theta _{s,t}(x)) - \Xi[    u       \cdot  \nabla   u     ](s,x)\Big ) ds	
	%	\nonumber \\
	%	&&
	+\int_0^t  \Xi[     u       \cdot  \nabla   u    ](s,x)	ds
	. \nonumber
\end{eqnarray*}
%In order to kill the contribution of the operator $\Xi$, and from the other hand to take $(\tau,\xi)=(t,x)$ for the regularisation of the remainder terms.
%As a consequence, we have to 
Let us choose $(\tau,\xi)=(t,x)$, next we apply $\P$, % to \eqref{Duhamel_u_FINAL}.
\begin{eqnarray*}%\label{ident_Pu}
	u     (t,x) &=&  \P [\hat P^{t,\cdot}_{0}  u_{0}    (t,\cdot)] (x) 
	+  \P [\hat G^{t,\cdot}_{0} \P f(t,\cdot)](x) 
	+  \P [ y \mapsto \hat  G^{t,y}_{0}  \big (    u_{\Delta}   [t,y]  \cdot  \nabla    u     \big )(t,y)]\Big |_{y=x}
	%\nonumber 
	\\
	&&+ \P [ y \mapsto  \hat  G^{t,y}_{0} \Big ( \Xi[    u       \cdot  \nabla   u    ](\cdot,\cdot)-\Xi[ u   \cdot  \nabla   u    ](\cdot,\theta _{\cdot,t}(y)) \Big )(t,y)	]\Big |_{y=x}
	\nonumber 
	\\
	&& 
	+ \P \Big [y \mapsto \int_0^t \Big ( \Xi[     u       \cdot  \nabla   u    ](s,\theta _{s,t}(y)) - \Xi[    u       \cdot  \nabla   u     ](s,y)\Big ) ds	
	\Big ] \Big |_{y=x}.
	%	\nonumber \\
	%	&&
	%\\
\end{eqnarray*}
From usual cancellation technique
%Comme par propriété de la fonction de Gauss,
\begin{equation*}
	\int_{\R^{d}} \hat  p^{t,x}  (s,t,x,y)  (y-u( \theta_{s,t }(x)))
	\ dy
	=0,
\end{equation*}
we eventually  derive the last Duhamel formula\footnote{Actually, this last argument is not essential for our result. It leads a better convergence speed in $n$, but requires more regularity on $u$.}
%De sorte qu'en injectant cette égalité dans \eqref{ident_Pu}
\begin{eqnarray}\label{ident_Pu_FINAL}
	u     (t,x) &=&  
	\P [\hat P^{t,\cdot}_{0}  u_{0}    (t,\cdot)] (x) 
	+  \P [\hat G^{t,\cdot}_{0} \P f(t,\cdot)](x) 
	+  \P [y \mapsto \hat  G^{t,y}_{0}\big (    u_{\Delta}   [t,y]  \cdot  \nabla    u     \big )(t,y)]\Big |_{y=x}
\nonumber 
	\\
	&&+ \P [ y \mapsto  \hat  G^{t,y}_{0} \Big ( \Xi[    u       \cdot  \nabla   u    ]^{\Delta^2}[t,\cdot]
	%](\cdot,\cdot)-\Xi[ u   \cdot  \nabla   u    ](\cdot,\theta _{\cdot,t}(y))
	\Big )(t,y)	]\Big |_{y=x}
	\nonumber 
	\\
	&& 
+ \P \Big [y \mapsto \int_0^t \Big ( \Xi[     u       \cdot  \nabla   u    ](s,\theta _{s,t}(y)) - \Xi[    u       \cdot  \nabla   u     ](s,y)\Big ) ds	
	\Big ] \Big |_{y=x},
\end{eqnarray}
with
\begin{equation*}%\label{def_Xi2}
	\Xi[    u       \cdot  \nabla   u    ]^{\Delta^2}[t,x](s,y)
	:=
	 \Xi[    u       \cdot  \nabla   u    ](s,y)
	 -\Xi[ u   \cdot  \nabla   u    ](s,\theta _{s,t}(x)) 
	 - (y-\theta _{s,t}(x)) \cdot \nabla \Xi[ u   \cdot  \nabla   u    ] (s,\theta _{s,t}(x))	.
\end{equation*}
Let us precise that $\nabla \Xi[ u   \cdot  \nabla   u    ] (s,\theta _{s,t}(x))$
%$\nabla \Xi[ u   \cdot  \nabla   u    ](t,x) $ 
stands for the gradient derivative of the function $x \mapsto \Xi[ u   \cdot  \nabla   u    ] (t,x)$ at the point $(s,\theta _{s,t}(x))$.

	\subsubsection{Time decomposition}
	\label{sec_decoup_temps}
	
	The goal, here, is to take advantage of the regularisation in small time intervals, like performed in \cite{hono_transport_v2:22}.
	For this purpose, we cut the time interval $[0,T]$ considered equation  \eqref{KOLMOLLI_xi} into ``small time pieces".

	Let $n \in \N$, and for each $k \in \leftB 0, n \rightB$, the time cutting is defined by
	\begin{equation*}
		t_k:= \frac{k}{n}T,
	\end{equation*}
	and the associated Cauchy problem for all  $x \in \R^d$,  $k \in \leftB 0, n-1 \rightB$ %$k \in \leftB 1, n-1 \rightB$ 
	and $t \in (t_{k}, t_{k+1}]$, 
	\begin{equation*}
%		\label{def_vn}
		\begin{cases}
			\partial_t  u_{k+1}(t,x)+ \P [u \cdot  \nabla   u_{k+1}](t,x) 
			%+ \frac nT \tilde G \Xi  ( u \cdot \nabla u )(t,x) \chi_{[t_{n-1},t_n]} (t)
			= \nu \Delta u_{k+1}(t,x)+ \mathbb P  f(t,x) 
			,\\
			u_{k+1}(t_k,x)= u_{k}(t_k,x),  %\, x \in \R^d .
		\end{cases}
	\end{equation*}
	with $u_{1}(0,\cdot)=u_0$, and for any $t \in (t_k, t_{k+1}]$, 
%	\begin{equation*}
		$u_{k+1}  (t,x):= u   (t,x)$.
		%\ \text{  et  } \omega_{k+1} (t,x):= \omega(t,x),
%	\end{equation*}
	From formula \eqref{ident_Pu_FINAL} \textit{mutatis mutandis}, for each $k \in \leftB 0,n-1 \rightB$
\begin{eqnarray}\label{ident_Pu_k1}
	u _{k+1}(t,x) &=&  	\P [\hat P_{t_{k}}^{t,\cdot}  u_{k}    (t_{k},\cdot)] (x) 
	+  \P [\hat G_{t_k}^{t,\cdot} \P f(t,\cdot)](x) 
	+  \P [ y \mapsto \hat  G_{t_k}^{t,y}  \big ( u_{\Delta} [t,y] \cdot  \nabla  u \big )(t,y)]\Big |_{y=x}
	\nonumber 
	\\
	&&+ \P [ y \mapsto  \hat  G_{t_k}^{t,y} \Big (
	\Xi[    u       \cdot  \nabla   u    ]^{\Delta^2}[t,\cdot]
	%\Xi[    u       \cdot  \nabla   u    ](\cdot,\cdot)-\Xi[ u   \cdot  \nabla   u    ](\cdot,\theta _{\cdot,t}(y)) 
	\Big )(t,y)	]\Big |_{y=x}
	%\nonumber 
	\\
&& + \P \Big [y \mapsto \int_{t_k}^t \Big ( \Xi[     u       \cdot  \nabla   u    ](s,\theta _{s,t}(y)) - \Xi[    u       \cdot  \nabla   u     ](s,y)\Big ) ds	
\Big ] \Big |_{y=x}
.	\nonumber 
	%\\
\end{eqnarray}

\subsubsection{Control of $u_{k+1}$ in  $ L^2$}
\label{sec_NS_L2}
%En utilisant l'estimée \eqref{ineq_reste1}
%%, \eqref{ineq_reste2}, \eqref{ineq_reste3} 
%dans la formule 
From triangular inequality and Duhamel formula \eqref{ident_Pu_k1},%,Duhamel_v_FINAL},
\begin{eqnarray*}%\label{ident_Pu_k1}
	\|	u _{k+1}(t,\cdot) \|_{L^2} &\leq &  	\| y\mapsto \P [\hat P_{t_{k}}^{t,\cdot}  u_{k}    (t_{k},\cdot)] (y)\|_{L^2} 
	%\P [	\hat P_{t_k}^{t,\cdot} u_k^{\Delta^2}(t_k,\cdot)[t,\cdot] (t,\cdot)+ \P [ u_k(t_k,\theta_{ t_k,t }(\cdot))] (x) 
	%	\nonumber \\
	%	&&
	+ \| y\mapsto  \P [\hat G_{t_k}^{t,y} \P f(t,\cdot)](y) \|_{L^2}
	\nonumber \\
	&&
	+ \| y\mapsto  \P [ \hat  G_{t_k}^{t,y}  \big ( u_{\Delta} [t,y] \cdot  \nabla  u \big )(t,\cdot)](y)\|_{L^2}
	%\nonumber 
	\\
	&&+ \Big \| y\mapsto  \P [ y \mapsto  \hat  G_{t_k}^{t,y} \Big (
	\Xi[    u       \cdot  \nabla   u    ]^{\Delta^2}[t,\cdot]
	%\Xi[    u       \cdot  \nabla   u    ](\cdot,\cdot)-\Xi[ u   \cdot  \nabla   u    ](\cdot,\theta _{\cdot,t}(y)) 
	\Big )(t,y)	]\Big |_{y=x} \Big \|_{L^2}
	\nonumber 
	\\
	&& + \Big \| x\mapsto  \P \Big [y \mapsto \int_{t_k}^t \Big ( \Xi[     u       \cdot  \nabla   u    ](s,\theta _{s,t}(y)) - \Xi[    u       \cdot  \nabla   u     ](s,y)\Big ) ds	
	\Big ] \Big |_{y=x} \Big \|_{L^2}.
	\nonumber 
	%\\
\end{eqnarray*}
From the Pythagorean theorem, we know that for any $\varphi \in L^2 (\R^d,\R^d)$,
\begin{equation*}
	\|\P \varphi \|_{L^2} \leq \|\varphi\|_{L^2}.
\end{equation*}
There is no deteriorating constant in the r.h.s. above, which is crucial for our analysis.
We provide a full comment on this argument in Section \ref{sec_comm_pyth} further, and why we cannot use a control in other $L^p$ space.
Hence,
\begin{eqnarray}\label{ineq_pythagore_uk}%\label{ident_Pu_k1}
\|	u _{k+1}(t,\cdot) \|_{L^2} &\leq &  	\| x\mapsto \hat P_{t_{k}}^{t,x}  u_{k}    (t_{k},\cdot) (x)\|_{L^2} 
%\P [	\hat P_{t_k}^{t,\cdot} u_k^{\Delta^2}(t_k,\cdot)[t,\cdot] (t,\cdot)+ \P [ u_k(t_k,\theta_{ t_k,t }(\cdot))] (x) 
%	\nonumber \\
%	&&
+ \| x\mapsto  \hat G_{t_k}^{t,x} \P f(t,x) \|_{L^2}
\nonumber \\
&&
+ \| x\mapsto   \hat  G_{t_k}^{t,x}  \big ( u_{\Delta} [t,x] \cdot  \nabla  u \big )(t,x)\|_{L^2}
%\nonumber 
%\\
%&&
+ \Big \| x\mapsto    \hat  G_{t_k}^{t,x} \Big (
\Xi[    u       \cdot  \nabla   u    ]^{\Delta^2}[t,\cdot]
%\Xi[    u       \cdot  \nabla   u    ](\cdot,\cdot)-\Xi[ u   \cdot  \nabla   u    ](\cdot,\theta _{\cdot,t}(y)) 
\Big )(t,x)	\Big \|_{L^2}
\nonumber 
\\
&& + \Big \| x\mapsto    \int_{t_k}^t \Big ( \Xi[     u       \cdot  \nabla   u    ](s,\theta _{s,t}(x)) - \Xi[    u       \cdot  \nabla   u     ](s,x)\Big ) ds	 \Big \|_{L^2}.
\nonumber 
%\\
\end{eqnarray}
Next, we have the important result of the change of variables,  
\begin{equation}\label{change_variable}
	\int_{\R^d} \varphi (\theta_{ s,t }(x)) dx = 	\int_{\R^d} \varphi (x) dx, 
\end{equation}
because the Jacobian determinant $(\det \nabla_x \theta_{ s,t }(x))$ % associated with divergence free
 matches with $1$, % identity matrix, 
 for more details see Lemma 1.1.1. in \cite{chem:98}.
In other words, we get
\begin{eqnarray}\label{ident_u_L2}
	&&
	\|	u _{k+1}(t,\cdot) \|_{L^2} 
	\nonumber \\
	&\leq &  	\| \tilde P_{t_{k-1}}  u_{k}    (t_{k},\cdot) \|_{L^2} 
	%\P [	\hat P_{t_k}^{t,\cdot} u_k^{\Delta^2}(t_k,\cdot)[t,\cdot] (t,\cdot)+ \P [ u_k(t_k,\theta_{ t_k,t }(\cdot))] (x) 
	%	\nonumber \\
	%	&&
	+ \|  \tilde G_{t_k} \P f(t,\cdot) \|_{L^2}
	\nonumber \\
	&&
	+ \Big \|x \mapsto \int_{t_k}^{t} \int_{\R^d} \tilde p(s,t,x,y)  \big ( (u(s,y)-u(s,x)) \cdot  \nabla  u (s,y)\big )  dy \, ds \Big \|_{L^2}
	\nonumber 
	\\
	&&+ \Big \| x \mapsto \int_{t_k}^{t} \int_{\R^d} \tilde p(s,t,x,y)  
	\Big (  	\Xi[    u       \cdot  \nabla   u    ] (s,y)
	-	\Xi[    u       \cdot  \nabla   u    ](s,x)
	 - (y-x) \cdot \nabla 	\Xi[    u       \cdot  \nabla   u   ](s,x) 
	 	\Big ) dy \, ds	\Big \|_{L^2}
%	 \\
%	 &&
%	\Xi[    u       \cdot  \nabla   u    ]^{\Delta^2}[t,\cdot]
%	%\Xi[    u       \cdot  \nabla   u    ](\cdot,\cdot)-\Xi[ u   \cdot  \nabla   u    ](\cdot,\theta _{\cdot,t}(y)) 
%	\Big )(t,\cdot) dy \, ds	\Big \|_{L^2}
	\nonumber 
	\\
	&& + 
	 \Big \| x\mapsto    \int_{t_k}^t \Big ( \Xi[     u       \cdot  \nabla   u    ](s,\theta _{s,t}(x)) - \Xi[    u       \cdot  \nabla   u     ](s,x)\Big ) ds	 \Big \|_{L^2}
	 %
%	\Big \| \tilde   G_{t_k}  \Big ( \Xi[     u       \cdot  \nabla   u    ](\cdot, \theta _{\cdot,t}(x)) - \Xi[    u       \cdot  \nabla   u     ](\cdot,y)\Big )(t,x) \Big \|_{L^2}.
	\nonumber \\
	&=:&
		\| \tilde P_{t_{k-1}}  u_{k}    (t_{k},\cdot) \|_{L^2}
			+ \|  \tilde G_{t_k} \P f(t,\cdot) \|_{L^2}
	+	\sum_{k=1}^3 \mathfrak{R}_k,
	%\\
\end{eqnarray}
recalling that the operators $\tilde P_{t_{k-1}}$ and $ \tilde G_{t_k}$ are defined in \eqref{def_tilde_G_tilde_P}, and where the remainder terms
by % are defined by
\begin{eqnarray}\label{def_Ri}
	\mathfrak{R}_1 &:=& \Big \|x \mapsto \int_{t_k}^{t} \int_{\R^d} \tilde p(s,t,x,y)  \big ( (u(s,y)-u(s,x)) \cdot  \nabla  u (s,y)\big )  dy \, ds \Big \|_{L^2},
	\nonumber \\
		\mathfrak{R}_2 &:=&
	\Big \| x \mapsto \int_{t_k}^{t} \int_{\R^d} \tilde p(s,t,x,y)  
	\Big (  	\Xi[    u       \cdot  \nabla   u    ] (s,y)
	-	\Xi[    u       \cdot  \nabla   u    ](s,x)
	- (y-x) \cdot \nabla 	\Xi[    u       \cdot  \nabla   u  ]  (s,x)
	\Big ) dy \, ds	\Big \|_{L^2},
		 \nonumber \\
		\mathfrak{R}_3 &:=&		 
			 \Big \| x\mapsto    \int_{t_k}^t \Big ( \Xi[     u       \cdot  \nabla   u    ](s,\theta _{s,t}(x)) - \Xi[    u       \cdot  \nabla   u     ](s,x)\Big ) ds	 \Big \|_{L^2}.
\end{eqnarray}
We insist that the change of variables is possible because the operator $\P$ is simplified thanks to the Pythagorean theorem.
\\

%\subsubsection
\textbf{Control of $\mathfrak{R}_1$}

By exponential absorbing property \eqref{ineq_absorb}, we readily derive
\begin{eqnarray*}%\label{ineq_R1}
	\mathfrak{R}_1 &\leq & \|\nabla u\|_{L^\infty_T L^\infty} \Big \|x \mapsto \int_{t_k}^{t} \int_{\R^d} \tilde p(s,t,x,y)  |y-x| \times |\nabla  u (s,y) |  dy \, ds \Big \|_{L^2}
	\nonumber \\
&\leq & C 
	\|\nabla u\|_{L^\infty_T L^\infty} \Big \|x \mapsto \int_{t_k}^{t} (t-s)^{\frac 12} \int_{\R^d} h_{C \nu } (t-s,x-y) |\nabla  u (s,y) |  dy \, ds \Big \|_{L^2},
%	\nonumber \\
%	&\leq & C 
%	\|\nabla u\|_{L^\infty_T L^\infty}  \int_{t_k}^{t} (t-s)^{\frac 12} \|\nabla  u (s,\cdot) \|_{L^2}  ds 
%	\nonumber \\
%	&\leq & C \frac{T^{\frac 32}}{n^{\frac 32}}
%	\|\nabla u\|_{L^\infty_T L^\infty} \|\nabla  u \|_{L^\infty_T L^2} 
%	\nonumber \\
%	&=:& \frac{T^{\frac 32}}{n^{\frac 32}} \mathbf N _{\eqref{ineq_R1}}(u).
\end{eqnarray*}
the heat kernel $ h_{C^{-1} \nu }$ is introduced in \eqref{def_h_tx}. We thus derive
\begin{equation}\label{ineq_R1}
	\mathfrak{R}_1
%	 &\leq & \|\nabla u\|_{L^\infty_T L^\infty} \Big \|x \mapsto \int_{t_k}^{t} \int_{\R^d} \tilde p(s,t,x,y)  |y-x| \times |\nabla  u (s,y) |  dy \, ds \Big \|_{L^2}
%	\nonumber \\
%	&\leq & C 
%	\|\nabla u\|_{L^\infty_T L^\infty} \Big \|x \mapsto \int_{t_k}^{t} (t-s)^{\frac 12} \int_{\R^d} h_{C^{-1} \nu } (t-s,x-y) |\nabla  u (s,y) |  dy \, ds \Big \|_{L^2}
%	\nonumber \\
	\leq  C 
	\|\nabla u\|_{L^\infty_T L^\infty}  \int_{t_k}^{t} (t-s)^{\frac 12} \|\nabla  u (s,\cdot) \|_{L^2}  ds 
	%\nonumber \\
	\leq  C \frac{T^{\frac 32}}{n^{\frac 32}}
	\|\nabla u\|_{L^\infty_T L^\infty} \|\nabla  u \|_{L^\infty_T L^2} 
%	\nonumber \\
	=: \frac{T^{\frac 32}}{n^{\frac 32}} \mathbf N _{\eqref{ineq_R1}}(u).
\end{equation}
We bear in mind that, for any $t \in [t_k,t_{k+1}]$, we have $t-t_k\leq t_{k+1}-t_k= \frac Tn$.
\\

%\subsubsection
\textbf{Control of $\mathfrak{R}_2$}

For this contribution, we need to separate the $L^2$ norm into two parts, one in $L^\infty$ which yields negligible coefficient in $n$, and another one for the space integration. 
That is to say, for  $\varepsilon \in (0,1)$, we get
\begin{eqnarray*}
%	&&\mathfrak{R}_2
%	\nonumber \\
	\mathfrak{R}_2 &\leq &
\Big \| x \mapsto \int_{t_k}^{t} \int_{\R^d} \tilde p(s,t,x,y)  
\Big (  	\Xi[    u       \cdot  \nabla   u    ] (s,y)
-	\Xi[    u       \cdot  \nabla   u    ](s,x)
- (y-x) \cdot \nabla 	\Xi[    u       \cdot  \nabla   u    ](s,x)
\Big ) dy \, ds	\Big \|_{L^\infty}^{\frac{\varepsilon}{2}}
\nonumber \\
&&\times
\Big \| x \mapsto \int_{t_k}^{t} \int_{\R^d} \tilde p(s,t,x,y)  
\Big (  	\Xi[    u       \cdot  \nabla   u    ] (s,y)
-	\Xi[    u       \cdot  \nabla   u    ](s,x)
- (y-x) \cdot \nabla 	\Xi[    u       \cdot  \nabla   u  ](s,x)
\Big ) dy \, ds	\Big \|_{L^{2-\varepsilon}}^{\frac{2-\varepsilon}{2}}
\nonumber \\
&\leq &
\| \nabla^2 	\Xi[    u       \cdot  \nabla   u   ] \|_{L^\infty_T L^\infty}^{\frac{\varepsilon}{2}}
\Big \| x \mapsto \int_{t_k}^{t} \int_{\R^d} \tilde p(s,t,x,y)  
|y-x|^2  dy \, ds	\Big \|_{L^\infty}^{\frac{\varepsilon}{2}}
\nonumber \\
&&\times 
\Big ( \int_{t_k}^{t} 2 \| 	\Xi[    u       \cdot  \nabla   u    ]  \|_{L^{2-\varepsilon}} + (t-s)^{\frac 12} \| \nabla \Xi[    u       \cdot  \nabla   u    ]  \|_{L^{2-\varepsilon}} ds \Big )^{\frac{2-\varepsilon}{2}}.
\end{eqnarray*}
Therefore, by Calderón–Zygmund inequality.
\begin{eqnarray}\label{ineq_R2}
%&&
	\mathfrak{R}_2
%	\nonumber \\
	&\leq &
%C 	\| \nabla^2 	\Xi[    u       \cdot  \nabla   u   ] \|_{L^\infty_T L^\infty}^{\frac{\varepsilon}{2}}
%	\Big (   \| 	\Xi[    u       \cdot  \nabla   u    ]  \|_{L^\infty_T L^{2-\varepsilon}} +\| \nabla \Xi[    u       \cdot  \nabla   u    ]  \|_{L^\infty_T L^{2-\varepsilon}}  \Big )^{\frac{2-\varepsilon}{2}}
%	\Big ( \int_{t_k}^t  (t-s) ds\Big ) ^{\frac{\varepsilon}{2}}
%	\nonumber \\
%	&&\times 
%	\Big ( \int_{t_k}^{t}  1+ (t-s)^{\frac 12} ds \Big )^{\frac{2-\varepsilon}{2}}
%	\nonumber \\
%		&\leq &
	C \frac{T^{1+\frac{\varepsilon}{2}}}{n^{1+\frac{\varepsilon}{2}}} 	\| \nabla^2 	\Xi[    u       \cdot  \nabla   u   ] \|_{L^\infty_T L^\infty}^{\frac{\varepsilon}{2}}
\nonumber \\
&&\times	\Big (   
	\| u \|_{L^\infty_T L^{\infty}}  
	\| \nabla   u    \|_{L^\infty_T L^{2-\varepsilon}}  
	+
	\| u \|_{L^\infty_T L^{\infty}}  
	\| \nabla^2   u    \|_{L^\infty_T L^{2-\varepsilon}}  
	+ 
	\|\nabla u \|_{L^\infty_T L^{\infty}}
	\| \nabla   u    \|_{L^\infty_T L^{2-\varepsilon}}   \Big )  ^{\frac{2-\varepsilon}{2}}
	\nonumber \\
	&=:&
	\frac{T^{1+\frac{\varepsilon}{2}}}{n^{1+\frac{\varepsilon}{2}}} 
	\mathbf N _{\eqref{ineq_R2}}^{(\varepsilon)}(u),
\end{eqnarray}

%\subsubsection
\textbf{Control of $\mathfrak{R}_3$}

Similarly, for any $\varepsilon \in (0,1)$, and from the definition of the flow $\theta _{s,t}$ in \eqref{def_theta},
\begin{eqnarray*}
		\mathfrak{R}_3 &\leq &		
		\Big \| x\mapsto    \int_{t_k}^t \Big ( \Xi[     u       \cdot  \nabla   u    ](s,\theta _{s,t}(x)) - \Xi[    u       \cdot  \nabla   u     ](s,x)\Big ) ds	  \Big \|_{L^\infty}^{\frac{\varepsilon}{2}}
		\nonumber \\
		&&\times
		\Big \| x\mapsto    \int_{t_k}^t \Big ( \Xi[     u       \cdot  \nabla   u    ](s,\theta _{s,t}(x)) - \Xi[    u       \cdot  \nabla   u     ](s,x)\Big ) ds	 	\Big \|_{L^{2-\varepsilon}}^{\frac{2-\varepsilon}{2}}
		\nonumber \\
		&\leq & 
		 	\| \nabla	\Xi[    u       \cdot  \nabla   u   ] \|_{L^\infty_T L^\infty}^{\frac{\varepsilon}{2}}
			\Big ( \int_{t_k}^t \int_{s}^{t} \|u(\tilde s ,\theta _{\tilde s,t}(\cdot))\|_{L^\infty} d\tilde s \, ds \Big )^{\frac{\varepsilon}{2}}
				\nonumber \\
			&&\times
			\Big (    \int_{t_k}^t  \| \Xi[     u       \cdot  \nabla   u    ](s,\theta _{s,t}(\cdot))\|_{L^{2-\varepsilon}} + \| \Xi[    u       \cdot  \nabla   u     ](s,\cdot)\|_{L^{2-\varepsilon}}  ds	 	\Big ) ^{\frac{2-\varepsilon}{2}}.
\end{eqnarray*}
From the change of variables principle for the divergence free flow, see \eqref{change_variable},
we deduce 
\begin{eqnarray}\label{ineq_R3}
	\mathfrak{R}_3 
	&\leq & 
2^{\frac{2-\varepsilon}{2}}	\| \nabla	\Xi[    u       \cdot  \nabla   u   ] \|_{L^\infty_T L^\infty}^{\frac{\varepsilon}{2}}
	\|u\|_{L^\infty_T L^\infty} ^{\frac{\varepsilon}{2}}
	\Big ( \int_{t_k}^t (t-s) ds \Big )^{\frac{\varepsilon}{2}}
%	\nonumber \\
%	&&\times
	\| \Xi[     u       \cdot  \nabla   u    ]\|_{L^\infty_T L^{2-\varepsilon}}  ^{\frac{2-\varepsilon}{2}}
	\Big (    \int_{t_k}^t  ds	 	\Big ) ^{\frac{2-\varepsilon}{2}}
	\nonumber \\
	&\leq &
C  \frac{T^{1+\frac{\varepsilon}{2}}}{n^{1+\frac{\varepsilon}{2}}}	\| \nabla	\Xi[    u       \cdot  \nabla   u   ] \|_{L^\infty_T L^\infty}^{\frac{\varepsilon}{2}}
	\|u\|_{L^\infty_T L^\infty} ^{\frac{\varepsilon}{2}}
	%	\nonumber \\
	%	&&\times
	\| u \|_{L^\infty_T L^{2-\varepsilon}}  ^{\frac{2-\varepsilon}{2}}
		\| \nabla   u    \|_{L^\infty_T L^{\infty}}  ^{\frac{2-\varepsilon}{2}}
%	 \nonumber \\
%	 &=:&
	 =:
	 \frac{T^{1+\frac{\varepsilon}{2}}}{n^{1+\frac{\varepsilon}{2}}} 
	 \mathbf N _{\eqref{ineq_R3}}^{(\varepsilon)}(u),
\end{eqnarray}
still by Calderón–Zygmund inequality.
\\

%\subsubsection
\textbf{Final control}

We readily obtain, from \eqref{ident_u_L2}, and by convolution principle,
\begin{equation}\label{ineq_u_k_L2_final}
	\|	u _{k+1}(t,\cdot) \|_{L^2} 
%\nonumber 
%\\
	\leq  	\|  u_{k}    (t_{k},\cdot) \|_{L^2} 
	%\P [	\hat P_{t_k}^{t,\cdot} u_k^{\Delta^2}(t_k,\cdot)[t,\cdot] (t,\cdot)+ \P [ u_k(t_k,\theta_{ t_k,t }(\cdot))] (x) 
	%	\nonumber \\
	%	&&
	+ \int_{t_k}^t \|   \P f(s,\cdot) \|_{L^2} ds
%	\nonumber \\
%	&&
	+ \frac{T^{\frac 32}}{n^{\frac 32}} \mathbf N _{\eqref{ineq_R1}}(u)
	+
\frac{T^{1+\frac{\varepsilon}{2}}}{n^{1+\frac{\varepsilon}{2}}} 
\big ( \mathbf N _{\eqref{ineq_R2}}^{(\varepsilon)}(u)+ \mathbf N _{\eqref{ineq_R3}}^{(\varepsilon)}(u)\big ).
%\nonumber
\end{equation}

\subsubsection{Control of $u$ in $L^\infty_T L^2$}

Iterating the above inequality gives
\begin{eqnarray*}%\label{ident_u_L2}
\|	u (T,\cdot) \|_{L^2}
% &=&	\|	u _{n}(t,\cdot) \|_{L^2} 
%\nonumber \\
\hspace{-0.25cm}		&\leq &  
		\hspace{-0.25cm}		\|  u_{1}    (t_{0},\cdot) \|_{L^2} 
			+
	\sum_{k=0}^{n-1} \Big [
	 \int_{t_k}^{t_{k+1}} \hspace{-0.25cm}	 \|   \P f(s,\cdot) \|_{L^2} ds
%	\nonumber \\
%	&&
	+ \frac{T^{\frac 32}}{n^{\frac 32}} \mathbf N _{\eqref{ineq_R1}}(u)
	+
	\frac{T^{1+\frac{\varepsilon}{2}}}{n^{1+\frac{\varepsilon}{2}}} 
	\big ( \mathbf N _{\eqref{ineq_R2}}^{(\varepsilon)}(u)+ \mathbf N _{\eqref{ineq_R3}}^{(\varepsilon)}(u)\big ) \Big ]
	\nonumber \\
	&=&
		\hspace{-0.25cm}		\|  u_{0} \|_{L^2} 
	+
 \int_{0}^t \|   \P f(s,\cdot) \|_{L^2} ds
%	\nonumber \\
%	&&
	+ \frac{T^{\frac 32}}{n^{\frac 12}} \mathbf N _{\eqref{ineq_R1}}(u)
	+
	\frac{T^{1+\frac{\varepsilon}{2}}}{n^{\frac{\varepsilon}{2}}} 
	\big ( \mathbf N _{\eqref{ineq_R2}}^{(\varepsilon)}(u)+ \mathbf N _{\eqref{ineq_R3}}^{(\varepsilon)}(u)\big ) 
	.
\end{eqnarray*}
Because $	\|	u (t,\cdot) \|_{L^2} 
$ does not depend on $n$, we can pass to the limit in the last identity.
% to derive the result.

%,  for any $t \in [0,T]$,
%\begin{equation*}%\label{ident_u_L2}
%	\|	u _{n}(t,\cdot) \|_{L^2} 
%\leq 
%	\|  u_{0} \|_{L^2} 
%	+
%	\int_{0}^t \|   \P f(s,\cdot) \|_{L^2} ds
%	.
%\end{equation*}

\section{Comments on Pythagorean argument}
\label{sec_comm_pyth}

Inequality \eqref{ineq_pythagore_uk} may be frustrating, we fail to adapt the analysis to the $L^p$ space, $p \neq 2$.
Indeed, by Calderón–Zygmund and from \eqref{ineq_u_k_L2_final}, we can expect an inequality of the type
%\begin{equation*}
%%&&	
%\|	u _{k+1}(t,\cdot) \|_{L^p} 
%%\nonumber \\
%	\leq  	\| x\mapsto \P [\hat P_{t_{k}}^{t,\cdot}  u_{k}    (t_{k},\cdot) (x)\|_{L^p} 
%	+ C \int_{t_k}^t \|   \P f(s,\cdot) \|_{L^p} ds
%%	\nonumber \\
%%	&&
%	+ \big (\frac{T^{\frac 32}}{n^{\frac 32}}
%	+
%		\frac{T^{1+\frac{\varepsilon}{2}}}{n^{1+\frac{\varepsilon}{2}}}  \big ) \mathbf N (u),
%\end{equation*}
%where $\mathbf N (u)>0$ is a constant depending on some norms of $u$. %Again by Calderón–Zygmund, we get,
\begin{eqnarray*}
	\|	u _{k+1}(t,\cdot) \|_{L^p} 
	&\leq &  C	\| \hat P_{t_{k}}^{t,\cdot}  u_{k}    (t_{k},\cdot) \|_{L^p} 
	+ C \int_{t_k}^t \|   \P f(s,\cdot) \|_{L^p} ds
	+ \big ( \frac{T^{\frac 32}}{n^{\frac 32}} 
	+
	\frac{T^{1+\frac{\varepsilon}{2}}}{n^{\frac{\varepsilon}{2}}} \big )  \mathbf N (u)
	\nonumber \\
		&\leq &  C	\|   u_{k}    (t_{k},\cdot) \|_{L^p} 
	+ C \int_{t_k}^t \|   \P f(s,\cdot) \|_{L^p} ds
	+ \big ( \frac{T^{\frac 32}}{n^{\frac 32}} 
	+
	\frac{T^{1+\frac{\varepsilon}{2}}}{n^{1+\frac{\varepsilon}{2}}}  \big )\mathbf N (u),
\end{eqnarray*}
where $\mathbf N (u)>0$ is a constant depending on some norms of $u$. 
%still by the change of variables associated with a flow of divergence free function.
The problem lies in the multiplicative constant in front of $\|   u_{k}    (t_{k},\cdot) \|_{L^p}$, which leads after iteration a geometric sequence, % depending on this constant,
\begin{eqnarray*}%\label{ident_u_L2}
	\|	u _{n}(t,\cdot) \|_{L^p} 
	&\leq &  
C^n	\|  u_{1}    (t_{0},\cdot) \|_{L^p} 
	+
	\sum_{k=1}^{n} C^k \Big [
C	\int_{t_k}^t \|   \P f(s,\cdot) \|_{L^p} ds
	+ \big ( \frac{T^{\frac 32}}{n^{\frac 32}} 
	+
	\frac{T^{1+\frac{\varepsilon}{2}}}{n^{1+\frac{\varepsilon}{2}}}  \big )\mathbf N (u)\Big ]
	.
\end{eqnarray*}
At the end of the day, there is no hope to make the remainder terms negligible and at the same time  the main terms convergent.

\section{Vorticity in $L^\infty_TL^1$}
\label{theo_vorticity}

%\subsection{Stat}
%\label{sec main}

%\subsection{Énoncé}

From now on, for a natural definition of the vorticity, we only consider the dimensional case $d=3$.
Adapting the above analysis, we derive a vorticity control first proved by \cite{cons:90}. 
\begin{THM}\label{THEO}
	
	%	\textcolor{red}{Dire que l'on suppose une unique soution rguliere, typiquement condition d'Oseen, ou petitesse de $u_0$.}
	%	
	Suppose that for a given  time $T>0$, there is a unique Leray solution $u$ of \eqref{Navier_Stokes_equation_v1} satisfying
	%	 is smooth enough, %	We mean by 	smooth when 
	%	specifically 
	$\|\nabla u\|_{L^\infty_T L^\infty} < + \infty$, then
	\begin{equation*}%\label{ineq_NS_THEO}
		\| \nabla \times u(T,\cdot)\|_{L^1}
		\leq 
		\|u_0\|_{L^1} +
		\int_0^T \| f(s,\cdot) \|_{L^1}  ds
		+	C	\|u_0\|_{L^2} +
		%\frac{1}{2}
		C	\int_0^T \| f(s,\cdot) \|_{L^2}  ds.
		% \Big ).
	\end{equation*}
	\end{THM}
	\begin{remark}
		The regularity condition on $u$ can be weaker, specifically if  $\sup_{t \in [0,T]}[u(t,\cdot)]_\gamma < + \infty$, for a given $\gamma \in (0,1]$, the conclusion is the same.
	\end{remark}
	%\begin{remark}
	%	contenu...
	%\end{remark}

	\subsection{Proof}
	%Démonstration : calculs \textit{a priori}}
\label{sec_preuve}

\subsubsection{Vorticity equation}

Let us apply the curl operator $\nabla \times$ to equation \eqref{Navier_Stokes_equation_v1},
\begin{equation}\label{Equation_omega}
	\partial_{t} \omega + u \cdot \nabla \omega = \nu \Delta \omega+ \omega \cdot  \nabla u + \nabla \times f,
\end{equation}
where the vorticity is determined by
%\begin{equation}%\label{def_omega}
$	\omega := \nabla \times u$ and $\omega_{0}=\nabla \times u_0$.
%\end{equation}
%called vorticity, for more details see for instance \cite{majd:bert:02}.
%\eqref{Navier_Stokes_equation_v2}
%
%The crucial advantage that we get from this equation, is the ``almost locality", unlike to \eqref{Navier_Stokes_equation_v1}, in a sens that there no Leray operator.
%This allows us to use the time decomposition trick introduced in \cite{hono_transport_v2:22}.

\subsubsection{Duhamel formula for $\omega$}

Thanks to the flow $\theta_{t,\tau} (\xi)$ defined in \eqref{def_theta}, we can rewrite equation \eqref{Equation_omega} by
\begin{eqnarray}
	\label{SMOLLI_xi}
	&&\partial_t  \omega  (t,x)
	+ 
	%\big [
	u (t,\theta_{t,\tau} (\xi)) 
	%+ (x-\theta_{ t,\tau } (\xi))\nabla u (t,\theta_{t,\tau} (\xi)) \big ]
	\cdot  \nabla  \omega(t,x) -   \nu \Delta  \omega (t,x)
	\nonumber \\
	&=&
	%	[  u  (t,\xi)-  u  ]
	%\P \big [
	u _{\Delta} [\tau,\xi](t,x) \cdot  \nabla  \omega
	%\big ]
	(t,x) 
	+ \omega \cdot  \nabla u(t,x)
	+ \nabla \times f(t,x),
	%\nonumber \\
\end{eqnarray}
where we recall that
%\begin{equation*}%\label{def_u_Delta}
	$u_{\Delta} [\tau,\xi](t,x) := 	   u(t,\theta_{t,\tau} (\xi))  - 	   u (t,x)$. 
%\end{equation*}
Furthermore, we can adapt Duhamel formula \eqref{Duhamel_u_INITIAL} to obtain,
\begin{equation*}%\label{Duhamel_u}
	\omega (t,x) =  \hat P_0^{\tau,\xi}  \omega_0(t,x) + \hat G_0^{\tau,\xi} \nabla \times f(t,x) 
	+  \hat  G_0^{\tau,\xi}  \big ( u _{\Delta} [\tau,\xi] \cdot  \nabla  \omega
	%\big ]
	+ \omega \cdot  \nabla u \big )(t,x),
\end{equation*}
Next, we can pick the \textit{freezing} point $(\tau,x) =(t,x)$,
\begin{equation}\label{Duhamel_omega}
	\omega (t,x) =  \hat P_0^{t,x}  \omega_0(t,x) + \hat G_0^{t,x} \nabla \times f(t,x) 
	+  \hat  G_0^{t,x}  \big ( u _{\Delta} [t,x] \cdot  \nabla  \omega
	%\big ]
	+ \omega \cdot  \nabla u \big )(t,x).
\end{equation}
We have to be careful, there is one contribution which is not purposed to be a remainder : $\omega \cdot  \nabla u $. The only way to overcome this difficulty is to use the regularity $L^2([0,T], H^1(\R^3, \R^3))$ of $u$.
We obtain the regularisation of the Gaussian kernel, namely we get 
\begin{eqnarray}\label{ineq_omega_Linfty_guide}
	\|\omega  (t,\cdot)\|_{L^1} &\leq& \|\omega_0\|_{L^1}  
	+ \int_0^t  \|\nabla \times  f(s,\cdot)\|_{L^1} ds 
	+C  \int_0^t  [\nu (t-s)]^{\frac 12}\| \nabla u \|_{L^\infty} \| \nabla  \omega(s,\cdot)\|_{L^1} ds
	\nonumber \\
	&&
	+  \int_0^t \| \omega \cdot  \nabla u(s,\cdot)\|_{L^1} ds
	\nonumber \\
	&\leq& \|\omega_0\|_{L^1}  
	+ \int_0^t \big ( \|\nabla \times  f(s,\cdot)\|_{L^1} 
	+ C  [\nu (t-s)]^{\frac 12}  \| \nabla u \|_{L^\infty} \| \nabla  \omega(s,\cdot)\|_{L^1} \big )ds
	\nonumber \\
	&&
	+  \frac 12 \int_0^t \| \omega(s,\cdot)|_{L^2} ^2+ \|  \nabla u(s,\cdot)\|_{L^2}^2 ds.
\end{eqnarray}
%où  $C_ {u }$ est une constante dépendant de $u $.
Let us pay attention that we can write $\|x \mapsto \hat P_0^{t,x}  \omega_0(t,x)  \|_{L^1} =\|\omega_0  \|_{L^1} $ by change of variables for a flow associated with an incompressible function, see \eqref{change_variable}.
%.  \cite{chem:98}. 
%\textcolor{blue}{A dvp pour plus de compréhension}

We cannot use a Gr\"onwall like lemma here, because  $\| \nabla u \|_{L^\infty}$ is not \textit{a priori} controlled.
In order to make negligible the unknown contribution $ \| \nabla u \|_{L^\infty} \| \nabla  \omega(s,\cdot)\|_{L^1} $, we use the time decomposition developed in Section \ref{sec_decoup_temps}.

\subsubsection{Time decomposition for the vorticity}

By the same notations introduced in Section \ref{sec_decoup_temps}, we can rewrite
the Cauchy problem \eqref{SMOLLI_xi} into small time intervals: for any $x \in \R^3$, and each $k \in \leftB 0, n-1 \rightB$,
\begin{equation*}%\label{Equation_omegak}
	\begin{cases}
		\partial_{t} \omega_{k+1}(t,x) + u \cdot \nabla \omega_{k+1} (t,x) = \nu \Delta \omega_{k+1}(t,x)+ \omega_{k+1} \cdot  \nabla u(t,x) + \nabla \times f(t,x),
		\ t \in (t_k,t_{k+1}],\\
		\omega_{k+1}(t_k,x)= \omega_k(t_k,x).
	\end{cases}
\end{equation*}
where
\begin{equation*}
	\omega_{k+1} (t,x):= \omega(t,x), \ 
	\omega_{0} (t,x):= \nabla \times u_0(x).
\end{equation*}
We can therefore modify Duhamel formula \eqref{Duhamel_omega} by
% into this time interval
%for any $t \in (t_{k}, t_{k+1}]$, by
\begin{equation*}%\label{Duhamel_omega_decomp}
	\omega_{k+1} (t,x) =  \hat P_{t_k}^{t,x}  \omega_k(t,x) + \hat G_{t_k}^{t,x} \nabla \times f(t,x) 
	+  \hat  G_{t_k}^{t,x}  \big ( u _{\Delta} [t,x] \cdot  \nabla  \omega_{k+1}
	%\big ]
	(t,x) 
	+ \omega_{k+1} \cdot  \nabla u \big )(t,x).
\end{equation*}
%and if $k=0$
%\begin{equation*}%\label{Duhamel_u}
%	\omega_{1} (t,x) =  \hat P_{0}^{t,x}  \omega_0(t,x) + \hat G_{0}^{t,x} \nabla \times f(t,x) 
%	+  \hat  G_{0}^{t,x}  \big ( u _{\Delta} [t,x] \cdot  \nabla  \omega_1
%	%\big ]
%	(t,x) 
%	+ \omega_1 \cdot  \nabla u \big )(t,x),
%\end{equation*}

\subsubsection{Control in $L^\infty_T L^1$ of $\omega$}

Adapting inequality \eqref{ineq_omega_Linfty_guide} to the interval $t \in (t_{k}, t_{k+1}]$, we obtain
\begin{equation*}%\label{ineq_omega_Linfty_guide}
	\|\omega_{k+1}  (t_{k+1},\cdot)\|_{L^1}
%	 &\leq& \|\omega_k\|_{L^1}  
%	+ \int_{t_k}^{t_{k+1}}  \|\nabla \times  f(s,\cdot)\|_{L^1} ds 
%	+C  \int_{t_k}^{t_{k+1}}  [\nu (t-s)]^{\frac 12}\| \nabla u \|_{L^\infty} \| \nabla  \omega(s,\cdot)\|_{L^1} ds
%	\nonumber \\
%	&&
%	+  \int_{t_k}^{t_{k+1}} \| \omega \cdot  \nabla u(s,\cdot)\|_{L^1} ds
%	\nonumber \\
%%	&\leq& \|\omega_k\|_{L^1}  
%%	+ \int_{t_k}^{t_{k+1}}  \big (\|\nabla \times  f(s,\cdot)\|_{L^1} 	+ C  [\nu (t_{k+1}-s)]^{\frac 12}  \| \nabla u \|_{L^\infty} \| \nabla  \omega(s,\cdot)\|_{L^1} \big ) ds
%	\nonumber \\
%	&&
%	+  \frac 12 \int_{t_k}^{t_{k+1}} \| \omega(s,\cdot)|_{L^2} ^2+ \|  \nabla u(s,\cdot)\|_{L^2}^2 ds
%	\nonumber \\
	\leq \|\omega_k\|_{L^1}  
	+ \int_{t_k}^{t_{k+1}} \hspace{-0.25cm}	 \big ( \|\nabla \times  f(s,\cdot)\|_{L^1}  +C\|  \nabla u(s,\cdot)\|_{L^2}^2
	+ C  [\nu (t_{k+1}-s)]^{\frac 12}  \| \nabla u \|_{L^\infty} \| \nabla  \omega(s,\cdot)\|_{L^1} \big ) ds,
%	\nonumber \\
%	&&
%	+ C \int_{t_k}^{t_{k+1}} \|  \nabla u(s,\cdot)\|_{L^2}^2 ds,
\end{equation*}
recalling that $\| \omega(s,\cdot)|_{L^2}\leq C \| \nabla u (s,\cdot)|_{L^2}$, by Calder\`on-Zygmund  inequality and Biot-Savart law.

Hence, to get the control at time $t=T$, we have to sum the previous inequality,
%il faut sommer l'inégalité précédente:
\begin{equation*}%\label{ineq_omega_Linfty_guide}
	%&&
	\|\omega  (T,\cdot)\|_{L^1}
	%\nonumber \\
%	& =& \|\omega _{n} (T,\cdot)\|_{L^1} 
%	\nonumber \\
	\leq
	%\|\nabla \times  \omega _{0}\|_{L^\infty} + \int_{t_k}^{t}  \|\nabla \times  f(s,\cdot)\|_{L^\infty} ds 
	%+ C_{u } \int_{t_k}^t [\nu(t-s)]^{\frac{1}{2}} \|\nabla   u  (s,\cdot)\|_{L^\infty} ds
	%\nonumber \\
	%&&+
	\|\nabla \times  \omega _{0}\|_{L^1} +
	\sum_{k=0}^{n-1}	%\Big \{ 
	\int_{t_k}^{t_{k+1}} \hspace{-0.25cm} \big (\|\nabla \times  f(s,\cdot)\|_{L^1} + C\|  \nabla u(s,\cdot)\|_{L^2}^2
	+ C  [\nu (t_{k+1}-s)]^{\frac 12}  \| \nabla u \|_{L^\infty} \| \nabla  \omega(s,\cdot)\|_{L^1} \big ) ds
%	\nonumber \\
%	&&
%	+ C\int_{t_k}^{t_{k+1}}  \|  \nabla u(s,\cdot)\|_{L^2}^2 ds \Big \}
	.
\end{equation*}
It is direct that for any  $s \in [t_k,t_{k+1}]$, $(t_{k+1}-s) \leq \frac T n$, and by Chasles relation,
\begin{eqnarray*}%\label{ineq_omega_Linfty_guide}
	\|\omega  (T,\cdot)\|_{L^1} 
	\leq
	%\|\nabla \times  \omega _{0}\|_{L^\infty} + \int_{t_k}^{t}  \|\nabla \times  f(s,\cdot)\|_{L^\infty} ds 
	%+ C_{u } \int_{t_k}^t [\nu(t-s)]^{\frac{1}{2}} \|\nabla   u  (s,\cdot)\|_{L^\infty} ds
	%\nonumber \\
	%&&+
	\|\nabla \times  \omega _{0}\|_{L^1} +
	\int_{0}^{T} \hspace{-0.25cm}	 \|\nabla \times  f(s,\cdot)\|_{L^1} ds 
+ C \|\nabla u \|_{L^2_TL^2}	+  \nu^{\frac{1}{2}} (\frac{T}{n})^{\frac 12} \int_{0}^{T}  \| \nabla u \|_{L^\infty} \| \nabla  \omega(s,\cdot)\|_{L^1}ds.
\end{eqnarray*}
The global solution $\omega  $ 
%and also $u  $ 
does not depend on  $n$, we can pass to the limit $n \to + \infty$,
% in the previous identity:
\begin{equation*}\label{ineq_omega_L1}
	\|\omega  (T,\cdot)\|_{L^1} 
	\leq
	\| \omega _{0}\|_{L^1} +
	\int_{0}^{T}  \|\nabla \times  f(s,\cdot)\|_{L^1} ds 
	+ C \|\nabla u \|_{L^2_TL^2}
%	\nonumber \\
%	&=:& \mathbf {N}_{\eqref{ineq_omega_L1}}(\omega_{0},f,T)
	.
\end{equation*}

\bibliographystyle{acm}
\bibliography{bibli}

\begin{thebibliography}{1}

\bibitem{chem:98}
{\sc Chemin, J.-Y.}
\newblock {\em Perfect incompressible fluids}, vol.~14 of {\em Oxford Lecture
  Series in Mathematics and its Applications}.
\newblock The Clarendon Press, Oxford University Press, New York, 1998.
\newblock Translated from the 1995 French original by Isabelle Gallagher and
  Dragos Iftimie.

\bibitem{cons:90}
{\sc Constantin, P.}
\newblock Navier-{S}tokes equations and area of interfaces.
\newblock {\em Comm. Math. Phys. 129}, 2 (1990), 241--266.

\bibitem{fuji:kato:64}
{\sc Fujita, H., and Kato, T.}
\newblock On the navier-stokes initial value problem. i.
\newblock {\em Archive for Rational Mechanics and Analysis 16\/} (1964),
  269--315.

\bibitem{giga:86}
{\sc Giga, Y.}
\newblock Solutions for semilinear parabolic equations in {$L^p$} and
  regularity of weak solutions of the {N}avier-{S}tokes system.
\newblock {\em J. Differential Equations 62}, 2 (1986), 186--212.

\bibitem{hono:NS_24}
{\sc Honor{\'e}, I.}
\newblock {Second derivatives of solutions to the 3D incompressible
  Navier-Stokes equation in Lebesgue spaces}.
\newblock working paper or preprint, Nov. 2024.

\bibitem{hono_transport_v2:22}
{\sc Honor{\'e}, I.}
\newblock {Transport equations in H{\"o}lder space by vanishing viscosity and
  applications}.
\newblock working paper or preprint.

\bibitem{lera:34}
{\sc Leray, J.}
\newblock Sur le mouvement d'un liquide visqueux emplissant l'espace.
\newblock {\em Acta Math. 63\/} (1934), 193--248.

\bibitem{osee:11}
{\sc {O}seen, C.~W.}
\newblock Sur les formules de {G}reen g\'{e}n\'{e}ralis\'{e}es qui se
  pr\'{e}sentent dans l'hydrodynamique et sur quelques-unes de leurs
  applications.
\newblock {\em Acta Math. 34}, 1 (1911), 205--284.

\bibitem{osee:12}
{\sc {O}seen, C.~W.}
\newblock Sur les formules de {G}reen g\'{e}n\'{e}ralis\'{e}es qui se
  pr\'{e}sentent dans l'hydrodynamique et sur quelques-unes de leurs
  applications.
\newblock {\em Acta Math. 35}, 1 (1912), 97--192.

\end{thebibliography}

\end{document}